\documentclass[a4paper]{amsart}
\usepackage{amsmath, amsfonts,amsthm,amssymb,amscd, verbatim,graphicx,color,multirow,booktabs, caption, mathdots,bm,chngcntr,adjustbox,longtable,mathtools,microtype}

\usepackage[top = 78pt, bottom = 72pt, inner=108pt, outer=108pt]{geometry}
\usepackage[dvipsnames, table]{xcolor}
\usepackage[backref]{hyperref}
\hypersetup{colorlinks=true,linkcolor=BurntOrange, citecolor=NavyBlue}

\usepackage[capitalise,noabbrev,nameinlink]{cleveref}
\usepackage[shortlabels]{enumitem}

\renewcommand{\wr}{ \mathop{\mathrm{wr}} }

\newtheorem{thm}{Theorem}[]

\theoremstyle{definition}

\newtheorem{remark}[thm]{Remark}

\newtheorem{problem}[thm]{Problem}

\title[Generators of groups acting arc-transitively]{On the number of generators of groups acting arc-transitively on graphs}

\author[M.~Barbieri]{Marco Barbieri}
\address{Dipartimento di Matematica ``Felice Casorati", University of Pavia, Via Ferrata 5, 27100 Pavia, Italy} 
\email{marco.barbieri07@universitadipavia.it}
\author[P.~Spiga]{Pablo Spiga}
\address{Dipartimento di Matematica e Applicazioni , University of Milano-Bicocca, Via Cozzi 55, 20125 Milano, Italy}
\email{pablo.spiga@unimib.it}

\thanks{The authors are members of the GNSAGA INdAM and PRIN ``Group theory and its applications'' research group and kindly acknowledge their support.}

\begin{document}	
	
	\begin{abstract}
		Given a finite connected graph $\Gamma$ and a group $G$ acting transitively on the vertices of $\Gamma$, we prove that the number of vertices of $\Gamma$ and the cardinality of $G$ are bounded above by a function depending only on the cardinality of $\Gamma$ and on the exponent of $G$. We also prove that the number of generators of a group $G$ acting transitively
		on the arcs of a finite graph $\Gamma$ cannot be bounded by a function of the valency alone.
	\end{abstract}
	
	\maketitle
	
	\section{Introduction}
	Let $\Gamma$ be a locally finite connected graph, let $G$ be a group acting transitively on the arcs of $\Gamma$. Suppose that, for any vertex $\alpha$ of $\Gamma$, the order of the vertex-stabilizer $G_\alpha$ is finite. In this note, an \emph{arc-transitive graph} is a pair $(\Gamma,G)$ with the properties just described. The \emph{local action of $(\Gamma,G)$} is the transitive action that a vertex-stabilizer $G_\alpha$ induces on the neighbourhood $\Gamma(\alpha)$ of the fixed point $\alpha$. In particular, the \emph{valency of $\Gamma$} coincides the degree of the local action of $(\Gamma,G)$.
	
	We say that a local action is \emph{graph-restrictive} if, for any arc-transitive graphs $(\Gamma,G)$ admitting such local action, then the order of the vertex-stabilizers is bounded from above by a constant depending only on the action. With this notation, the famous Weiss Conjecture (posed in \cite{Weiss1978}) states that primitive groups are graph-restrictive. We refer to \cite{PotocnikSpigaVerret2012} for an extensive study of this notion.
	
	Let $(\Gamma,G)$ be an arc-transitive graph of valency $d$ whose local action is graph-restrictive. One can easily check that, by a connectedness argument (see the proof of Thereom~\ref{thm:main}), the group $H$ generated by those elements $d$ of $G$ sending a prescribed vertex $\alpha$ to all its neighbours in turn is transitive. Therefore, by Frattini's Argument, $G=G_\alpha H$. In particular, $|G_\alpha|+d$ elements are sufficient to generate $G$. Hence the minimal number of generators for $G$, denoted by $\mathbf{d}(G)$ can be bounded with a function of the local action alone.
	
	More surprisingly, for every arc-transitive graph $(\Gamma,G)$ of valency at most $4$, the minimal number of generators of $G$ is bounded by a constant regardless of the local action. The result is trivial for $d\in\{1,2\}$. For $d\in\{3,4\}$, some deeper concepts enter the picture. For every arc-transitive graph $(\Gamma,G)$, there is a universal cover of the form $(\mathcal{T}_d, G_\alpha \ast_{G_{\alpha\beta}} G_{\{\alpha,\beta\}})$, where $\alpha$ and $\beta$ are two adjacent vertices of $\Gamma$, adn $\mathcal{T}_d$ is the infinite tree of valency $d$. For every amalgamated product appearing in the universal cover for valency $3$ and $4$, an explicit presentation has been produced: see \cite{Djokovic,DjokovicMiller,Potocnik2009}. A direct inspection of these presentations proves that the minimal number of generators for $G$ is at most $10$.
	
	One could dare to conjecture that there exists a function $\mathbf{f}: \mathbb{N} \to \mathbb{N}$ that takes the valency of the graph $\Gamma$ as input, and returns an upper bound for $\mathbf{d}(G)$.
	The main contribution of this note is proving that such a function cannot exist.
	
	\begin{thm}\label{thm:counter}
		There exists no function $\mathbf{f}: \mathbb{N} \to \mathbb{N}$ such that, for every arc-transitive graph $(\Gamma,G)$ of valency $d$,
		\[\mathbf{d}(G) \le \mathbf{f}(d) \,.\]
	\end{thm}

	\begin{remark}
		To prove the veracity of Theorem~\ref{thm:counter}, we will exhibit an infinite family $\mathcal{F}$ of pairs $(\Gamma_h,G_h)$ such that the valency of the graphs is a constant (at least $8$), while $\mathbf{d}(G_h)$ grows linearly with the exponent of the group. We would like to remark that the graphs $\Gamma_h$ carry an oustanding similiarity with those build in \cite{HujdurovicPotocnikVerret,PotocnikSpigaVerret2014,PotocnikSpigaVerretSpectrum} to prove that, for some imprimitive local actions of degree $6$, the order of the vertex-stabilizers grows exponentially.
	\end{remark}

	We also observe that, in our construction, $G$ is not the automorphism group of $\Gamma$. This prompts the following question.
	\begin{problem}\label{problem:Aut}
		Is there a function $\mathbf{f}: \mathbb{N} \to \mathbb{N}$ such that, if $\Gamma$ is a connected arc-transitive graph of valency $d$, then
		\[ \mathbf{d}\left( \mathrm{Aut}(\Gamma) \right) \le \mathbf{f}(d) \,?\]
	\end{problem}

	Moreover, for our current and limited knowledge of arc-transitive graphs $(\Gamma, G)$, having $\mathbf{d}(G)$ bounded appears to be quite common. Therefore, we ask the following.
	\begin{problem}
		Which assumption on the arc-transitive graph $(\Gamma,G)$ are necessary to bound $\mathbf{d}(G)$ with a function of the valency (or the local action)?
	\end{problem}
	
	To conclude, we give a bound on the order of the group $G$ appearing in an arc-transitive graph $(\Gamma,G)$ depending on the valency $d$ and the exponent of $G$. Our proof does not rely on the fact that the local action is transitive, thus the result holds for the larger class of vertex-transitive graphs. In this note, a \emph{vertex-transitive graph} is a pair $(\Gamma,G)$ where $\Gamma$ is a finite connected graph and $G$ is a subgroup of $\mathrm{Aut}(\Gamma$) whose action on the vertex-set of $\Gamma$ is transitive.
	
	We underline that the exponent of $G$ is not a local feature of the graph. For instance, if $\Gamma$ is a cycle of odd length, then the local action is isomorphic to $C_2$, while the exponent of $\mathrm{Aut}(\Gamma)$ is twice the length of the cycle. Finally, we point out that the function $\mathbf{g}$ appearing in Theorem~\ref{thm:main} is not explicitly determined, but it depends on the solution of the Burnside Restricted Problem. Moreover, even for the few cases in which we can make it explicit, its growth is at least factorial in both parameters.
	
	\begin{thm}\label{thm:main}
		Let $(\Gamma,G)$ be a vertex-transitive graph. Suppose that the valency of $\Gamma$ is $d$, and that the exponent of $G$ is $e$. Then there exists a function $\mathbf{g}: \mathbb{N}\times \mathbb{N} \rightarrow \mathbb{N}$ such that
		\[|G|\leq \mathbf{g}(d,e) \,.\]
	\end{thm}

	\section{Proof of Theorem~\ref{thm:counter}}
	
	Let $h$ be a positive integer, and let $p$ be a prime. Set
	\[H=C_{p^h}\times C_{p^h}= \langle a,b \mid a^{p^h} = b^{p^h} = [a,b] = 1 \rangle \,,\]
	and consider the group algebra $\mathbb{F}_p[H]$. Define recursively the following chain of $H$-modules 
	\[ \gamma_0 := \mathbb{F}_p[H], \quad \hbox{and, \, for any } i\ge 1, \]
	\[ \gamma_i := [\gamma_{i-1}, H] = \langle v-vh \mid v \in \gamma_{i-1}, h\in H \rangle_{\mathbb{F}_p}. \]
	Consider the $(p^h \times p^h)$-square unipotent matrix  
	\[M=\begin{pmatrix}
		1 & 1 & 0 & \dots & 0 & 0 & 0 \\
		0 & 1 & 1 & \dots & 0 & 0 & 0 \\
		0 & 0 & 1 & \dots & 0 & 0 & 0 \\
		\vdots & \vdots & \vdots & \ddots &  \vdots & \vdots & \vdots \\
		0 & 0 & 0 & \dots & 1 & 1 & 0 \\
		0 & 0 & 0 & \dots & 0 & 1 & 1 \\
		0 & 0 & 0 & \dots & 0 & 0 & 1 \\
	\end{pmatrix},\]
	and denote the identity matrix of size $p^h$ by $id$. Observe that $M$ is an element of order $p^h$ in $\mathrm{GL}_{p^h}(p)$.
	
	The natural basis for the group algebra $\mathbb{F}_p[H]$ is
	\[ \left(a^ib^j \mid i,j \in \left\{0,1,\ldots,p^h-1\right\} \right) \,.\]
	We have that
	\[ \mathcal{B} = \left(e_{xy}=(a-1)^x(b-1)^y \mid x,y \in \left\{0, 1, \ldots, p^h-1\right\} \right)\]
	is also a basis. For convenience, we set $e_{xp^h} = e_{p^hy} = 0$, for every $x,y \in \left\{0, 1, \ldots, p^h-1\right\}$.  Observe that
	\[ \begin{split}
		e_{xy} \cdot a &= (a-1)^x(b-1)^y \cdot a \\&= (a-1)^x(1+a-1)(b-1)^y
		\\&= (a-1)^x(b-1)^y + (a-1)^{x+1}(b-1)^y
		\\&= e_{xy} + e_{(x+1)y} \,,
	\end{split}\]
	and
	\[ \begin{split}
		e_{xy} \cdot b &= (a-1)^x(b-1)^y \cdot b \\&= (a-1)^x(b-1)^y(1+b-1)
		\\&= (a-1)^x(b-1)^y + (a-1)^x(b-1)^{y+1}
		\\&= e_{xy} + e_{x(y+1)} \,.
	\end{split}\]
	Therefore, in the basis $\mathcal{B}$, the action of $a$ and $b$ on $\mathrm{F}_p[H]$ is given by the matrices by the matrices $ M \otimes id$ and $ id \otimes M$ respectively.
	
	
	With a direct computation, we get that
	\[\gamma_i = \langle e_{xy} \mid x+y \ge i \rangle_{\mathbb{F}_p}, \quad \hbox{and}\]
	\[\gamma_i / \gamma_{i+1} = \langle e_{xy} + \gamma_{i+1} \mid x+y = i \rangle_{\mathbb{F}_p}.\]
	Recall that, for any $H$-module $V$, we denote by $\mathbf{d}_H(V)$ the minimal number of generators of $V$ as an $H$-module. Since, by construction, $\gamma_i / \gamma_{i+1}$ is a trivial section of $\mathbb{F}_p[H]$, we have that
	\[ \mathbf{d}_H \left( \gamma_i / \gamma_{i+1} \right) = \mathrm{dim}_{\mathbb{F}_p} \left( \gamma_i / \gamma_{i+1} \right) = \min \{i+1, 2p^h-i-1 \}.\]
	Therefore, $\mathbf{d}_H(\gamma_{p^h-1}) = p^h$.
	In particular, 
	\[\mathbf{d} \left(\gamma_{p^h-1}\rtimes H\right) = p^h+2 \,.\]
	
	The group $H$ acts regularly on the vertex-set of the $4$-valent Cayley graph defined by
	\[\Delta:=\mathrm{Cay}(H,\{a,a^{-1},b,b^{-1}\}).\]
	Recall that, for any two graphs $\Gamma,\Delta$, the \emph{wreath product of $\Delta$ by $\Gamma$}, denoted by $\Gamma \wr \Delta$, is the graph having vertex-set $V\Gamma \times V\Delta$, where $(\delta_1,\gamma_1)$ and $(\delta_2,\gamma_2)$ are adjacent if either $\delta_1 = \delta_2$ and $\{\gamma_1,\gamma_2\}$ is an edge of $\Gamma$, or $\{\delta_1, \delta_2\}$ is an edge of $\Delta$. We can finally define $\Gamma_h$ as the wreath product of the empty graph on $p$ vertices, $p\mathbf{K}_1$, by the Cayley graph $\Delta$, that is, $\Gamma_h := p\mathbf{K}_1 \wr \Delta$. Note that, unless $p=2$ and $h=1$, $\Gamma_h$ has valency $4p$.
	
	
	Observe that $\gamma_{p^h-1}\rtimes H$ embeds into $C_p \mathrm{wr} H$, which is a subgroup of the automorphism group of $\Gamma_h$. Moreover, as $\gamma_{p^h-1}$ is not trivial, $\gamma_{p^h-1}\rtimes H$ must be transitive on the vertices of $\Gamma_h$. On the other hand, since $\gamma_{p^h-1}\rtimes H$ preserves the lifting of the labels $\{a,a^{-1},b,b^{-1}\}$ from the Cayley graph $\Delta$, this action is not arc-transitive. In particular, the local action of $(\Gamma_h, \gamma_{p^h-1}\rtimes H)$ is intransitive with $4$ distinct orbits.
	
	We need to extend $H$ with some outer automorphisms to achieve arc-transitivity. Consider the automorphism of $H$ defined on the generators by
	\[ \varphi:\, a \mapsto b,\, b\mapsto a, \quad \hbox{and} \quad \psi:\, a\mapsto a^{-1},\, b\mapsto b^{-1} \,. \]
	Observe that both are involution, and that $\langle \varphi,\psi \rangle$ is isomorphic to the Klein group. Extend the multiplication on $\mathbb{F}_p[H]$ by putting, for every $\varepsilon,\delta \in \mathbb{Z}_2$,
	\[ \left(\sum_{h\in H} \lambda_h h \right) \cdot \left(\varphi^\varepsilon \psi^\delta\right) = \sum_{h\in H} \lambda_h h ^ {\varphi^\varepsilon \psi^\delta} \,. \]
	With this operation, $\gamma_0 = \mathbb{F}_p[H]$ is an $\mathbb{F}_p[H\rtimes\langle\varphi,\psi\rangle]$-module.
	Furthermore, for all indices $i$, since $\gamma_i$ is a characteristic subgroup of $\gamma_0 \rtimes H$, then $\gamma_i$ is an $\mathbb{F}_p[H\rtimes\langle\varphi,\psi\rangle]$-submodule. In particular, the group $G_h := \gamma_{p^h-1} \rtimes (H\rtimes \langle \varphi, \psi \rangle )$ is a subgroup of $\mathrm{Aut}(\Gamma_h)$. Finally, as the index of $H$ in $H\rtimes\langle\varphi,\psi\rangle$ is $4$, we have that
	\[\mathbf{d}_{H\rtimes\langle\varphi,\psi\rangle} (\gamma_{p^h-1}) \ge \frac{p^h}{4} \,.\]
	Hence, we obtain that
	\[ \mathbf{d}(G_h) \ge \frac{p^h}{4} + 3 \,.\]
	
	Let us go back to the Cayley graph $\Delta$. Observe that $\langle \varphi, \psi \rangle$ is transitive on the connection set of $\Delta$. This implies that $H \rtimes \langle \varphi,\psi \rangle$ is an arc-transitive subgroup of $\mathrm{Aut}(\Delta)$. Therefore, the local action of $(\Gamma_h,G_h)$ transitively permutes the $4$ orbits defined by $(\Gamma_h,\gamma_{p^h-1}\rtimes H)$. Hence, the action of $G_h$ on the graph $\Gamma_h$ is arc-transitive.
	
	Therefore, for a fixed prime $p$, the family
	\[ \mathcal{F}_p = \{ (\Gamma_h,G_h) \mid h\ge 2 \}\]
	contains only graphs of valency $4p$, meanwhile
	\[\lim\limits_{h \rightarrow +\infty} \mathbf{d}(G_h) \ge \lim\limits_{h \rightarrow +\infty} \frac{p^h}{4} + 3 = + \infty.\]
	This concludes the proof of Theorem~\ref{thm:counter}.

	\section{Proof of Theorem~\ref{thm:main}}
	
	Let $(\Gamma,G)$ be a vertex-transitive graph, and suppose that $d$ is the valency of $\Gamma$.  Choose a vertex $\alpha$ of $\Gamma$, and, for any of its neighbours $\beta$, let $g_\beta\in G$ be an element such that
	\[\alpha^{g_\beta}=\beta \,.\]
	(This element exists by the transitivity of $G$ on the vertex-set of $\Gamma$.)
	Define the subgroup
	\[ H:= \langle g_\beta \mid \beta \hbox{ is a neighbour of } \alpha \rangle \leq G \,.\]
	By a routine connectedness argument, $H$ is transitive on the vertex-set of $\Gamma$. Hence, by Frattini's Argument, $G=G_\alpha H$.
	
	Let $\mathbf{B}(d,e)$ be the function solving the Burnside Restricted Problem for a finite group with $d$ generators and exponent $e$. The existence of this function for all the choices of $d$ and $e$ was proved by Zel'manov in \cite{ZelmanovOdd,ZelmanovEven}. Then, $|H|\leq \mathrm{B}(d,e)$.
	
	We now need to bound the order of $G_\alpha$.  Note that $G_\alpha$ can be embedded into $\mathrm{Sym}(V\Gamma  -  \{\alpha\})$, hence, as $|V\Gamma|\leq |H|$,
	\[ |G_\alpha| \leq (\mathbf{B}(d,e)-1)! \,.\]
	Therefore, we obtain that
	\[ |G| = |H||G_\alpha| \leq \mathbf{B}(d,e)(\mathbf{B}(d,e)-1)! \,. \]
	
	\bibliographystyle{plain}
	\bibliography{bibGenerators}
\end{document}